\documentclass[11pt]{amsart}
\pdfoutput=1
\usepackage{graphicx}
\usepackage{wasysym}
\usepackage{amssymb, amsmath, amscd}
\usepackage{float}
\usepackage[mathscr]{euscript}
\usepackage{array}
\usepackage{chngcntr}
\usepackage[bookmarks=false]{hyperref}

\newtheorem{theorem}{Theorem}[section]

\newtheorem{lemma}[theorem]{Lemma}

\newtheorem{proposition}[theorem]{Proposition}

\newtheorem{definition-proposition}[theorem]{Definition-Proposition}
\newtheorem{lemma-notation}[theorem]{Lemma-Notation}

\theoremstyle{definition}
\newtheorem{definition}[theorem]{Definition}

\newtheorem{remark}[theorem]{Remark}

\newcommand{\Z}{\mathbb{Z}}

\newcommand{\C}{\mathbb{C}}

\newcommand{\twopartdef}[4]
{
	\left\{
		\begin{array}{ll}
			#1 & \mbox{if } #2 \\
			#3 & \mbox{if } #4
		\end{array}
	\right.
}

\newcommand{\twopartdefotherwise}[3]
{
	\left\{
		\begin{array}{ll}
			#1 & \mbox{if } #2 \\
			#3 & \mbox{\textrm{otherwise.}}
		\end{array}
	\right.
}

\begin{document}
\title{Hurwitz theory of elliptic orbifolds, II}
\author{Philip Engel}
\address[Philip Engel]{University of Georgia}
\email{philip.engel@uga.edu}
\maketitle

\begin{abstract} An {\it elliptic orbifold} is the quotient of an elliptic curve by a finite group. In 2001, Eskin and Okounkov proved that generating functions for the number of branched covers of an elliptic curve with specified ramification are quasimodular forms for $SL_2(\Z)$. In 2006, they generalized this theorem to the enumeration of branched covers of the quotient of an elliptic curve by $\pm 1$, proving quasi-modularity for $\Gamma_1(2)$. In 2017, the author generalized their work to the quotient of an elliptic curve by $\langle \zeta_N\rangle$ for $N=3, 4, 6$, proving quasimodularity for $\Gamma_1(N)$.

In these works, both Eskin-Okounkov and the author had to assume that there was at least one orbifold point of order $N$ over which there was no ramification. Here we remove that assumption, with the caveat that the generating functions are only quasimodular for $\Gamma(N)$. We deduce the following corollary: Let $h_6(\vec{\kappa},q)$ be the generating function whose $q^n$ coefficient is the number of surface triangulations with $2n$ triangles, such that the set of non-zero curvatures is $\kappa_i$. Here the {\it curvature} of a vertex is six minus its valence. Then under the substitution $q=e^{2\pi i \tau}$, the function $h_6(\vec{\kappa},q)$ is a quasimodular form for $\Gamma_1(6)$ with weight bounded in terms of $\vec{\kappa}$. This statement in turn implies that the Masur-Veech volume of any stratum of sextic differentials is polynomial in $\pi$. \end{abstract}

\section*{Introductory remarks}

This work is a continuation of \cite{engel} from 2017, which proved a special case of the main theorem here. Rather than significantly revise \cite{engel}, I decided it would be more useful to keep the precursor to this paper available, rather than to hide its origins. In addition, it provides much more background. As described in \cite{engel}, the main application of Theorem \ref{modularity} is that the generating function for flat surfaces tiled by a unit equilateral triangle in a stratum of sextic differentials lies in a ring of quasi-modular forms for $\Gamma_1(6)$. I recommend starting with \cite{engel} first, if flat surfaces are your field of interest. On the other hand, for those who are unconcerned with the enumeration of tilings, this paper stands alone as the computation of the $q$-trace of certain operators on the fermionic Fock space.

More precisely, we compute the trace of a product of vertex operators on $\mathbb{L}^{\otimes N}$, where $\mathbb{L}$ is the {\it charge zero subspace of Fock space}---a vector space with a basis $v_\lambda$ indexed by partitions $\lambda=\{\lambda_1\geq\lambda_2\geq\dots\}$. Define the operators $H$ and $\mathcal{E}_0(z)$ on $\mathbb{L}$ by the actions \begin{align*} Hv_\lambda &= |\lambda|v_\lambda \\ \mathcal{E}_0(z)v_\lambda &= \sum_i e^{\lambda_i-i+\frac{1}{2}}v_\lambda.\end{align*} Then, the trace that we compute is $$tr_{\mathbb{L}^{\otimes N}}\!\!\bigotimes_{r\textrm{ mod }N} \!\!\!\left(q^H\mathcal{E}_0(z_{r,1})\dots\mathcal{E}_0(z_{r,{n_r}})\Psi_r\right) \circ \mathfrak{R}_N $$ where $\Psi_r$ is a vertex operator with specified coefficients and $\mathfrak{R}_N$ is the operator which cyclically rotates the tensor factors of $\mathbb{L}^{\otimes N}$. The resulting trace is a contour integral of a function expressed in terms of Jacobi theta functions of level $\Gamma(N)$ in the $z$ variables and $q$. We conclude that the Taylor coefficients in the $z_{r,i}$ variables are quasimodular forms. Quasimodularity of generating functions of triangulations follows from the $N=6$ case, generalizing the result in \cite{engsmi} to all genera.

\section{Vertex operators associated to quotients and cores}

\begin{figure}
\makebox[\textwidth][c]{\includegraphics[width=5.5in]{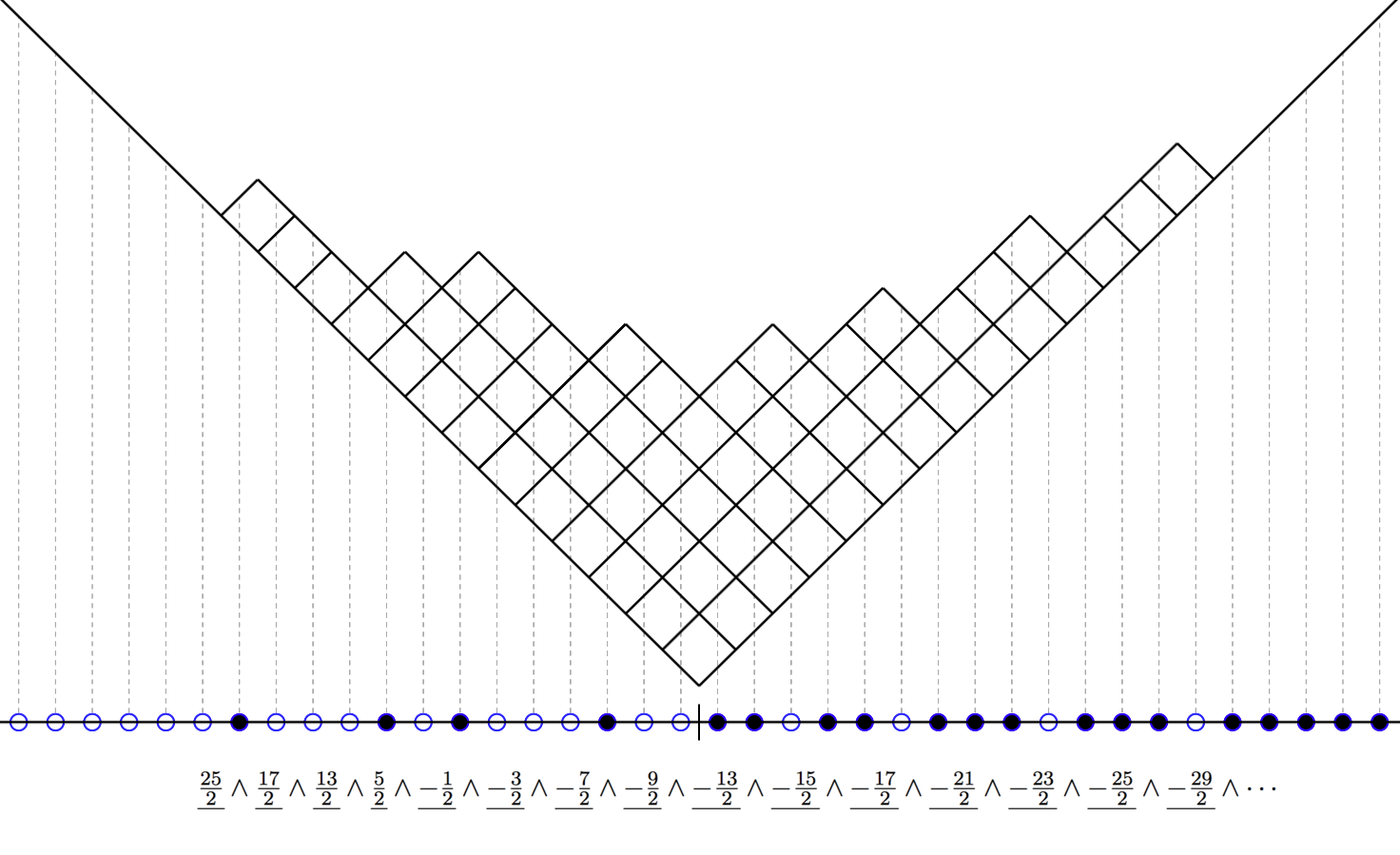}}
\caption{The Young diagram of a partition $\lambda$, transformed into the associated half-infinite wedge $v_\lambda$. We think of $v_\lambda$ as an infinite sequence of beads on a string, with possible positions at the set of half-integers $\frac{1}{2}+\Z$, ordered in decreasing order left to right. All slots sufficiently far to the left are empty, and all slots sufficiently far to the right are occupied.}
\label{fig1}
\end{figure}

Recalling the notations in Sections 3, 4, and 5 of \cite{engel}. Denote by $V$ a vector space with a basis indexed by half-integers \begin{align*}V:=\textrm{span}\,\{\underline{i}\,:\,i\in\tfrac{1}{2}+\Z\}.\end{align*} The {\it half-infinite wedge} or {\it Fock space} $\Lambda^{\infty/2}V$ is the vector space spanned by all formal symbols \begin{align*}\underline{i_1}\wedge\underline{i_2}\wedge \underline{i_3}\wedge\cdots\end{align*} such that $i_j>i_{j+1}$ and $i_{j+1}=i_j-1$ for all $j\gg 0$. See \cite{rz2} for an introduction to the half-infinite wedge. Define an inner product  $\langle \cdot \big{|}\cdot \rangle$ on Fock space by declaring these symbols orthonormal. We call this basis the {\it fermionic basis}. The {\it charge $C$ subspace} is the subspan of symbols for which $i_j=-j+\frac{1}{2}+C$ for all $j\gg 0$. The {\it shift operator} $S$ on Fock space increases all $i_j$ by $1$. Then $S$ sends the charge $C$ subspace isomorphically to the charge $C+1$ subspace. The charge zero subspace $\mathbb{L}$ has a basis indexed by partitions $\lambda=\{\lambda_1\geq \lambda_2\geq \lambda_3\geq \dots\}$ of all integers: \begin{align*}v_\lambda:=\underline{\lambda_1-1+\tfrac{1}{2}}\wedge\underline{\lambda_2-2+\tfrac{1}{2}}\wedge\underline{\lambda_3-3+\tfrac{1}{2}}\wedge\cdots.\end{align*} See Figure \ref{fig1}. There is a tensor product decomposition by charge \begin{align*}\Lambda^{\infty/2}V=\mathbb{L}\otimes \C[\Z]\end{align*} such that the fermionic basis can be written as $v_\lambda\otimes C$, the charge $C$ subspace is $\mathbb{L}\otimes C$, and the shift acts by $C\mapsto C+1$ on the second tensor factor.

Define the {\it energy operator} $H$ on Fock space by the action \begin{align*} H(v_\lambda\otimes C) = (|\lambda|+\tfrac{1}{2}C^2)\,v_\lambda\otimes C.\end{align*} Then the eigenspace of $H\big{|}_\mathbb{L}$ with eigenvalue $n$ is spanned by $v_\lambda$ as $\lambda$ ranges over the partitions of $n$. 
%This action extends to all of Fock space by setting \begin{align}H(v_\lambda\otimes C)= \left(|\lambda| + \tfrac{1}{2}C^2\right)v_\lambda\otimes C.\end{align} 
Recall that the {\it Heisenberg algebra} is the Lie algebra \begin{align*}\C e\oplus \bigoplus_{n\neq 0} \C\alpha_n\end{align*} such that $e$ is central and $[\alpha_n,\alpha_m]=n\delta_{n,-m}e$. Define an action of the Heisenberg algebra on $V$ by $\alpha_n(\underline{i})=\underline{i-n}.$ There is an induced Lie algebra action on Fock space by   \begin{align}\begin{aligned}\label{action} \alpha_n\,:\,\underline{i_1}\wedge \underline{i_2}\wedge \underline{i_3}\wedge\dots\mapsto &\, \underline{i_1-n}\wedge \underline{i_2}\wedge \underline{i_3}\wedge\dots+ \\ & \,\underline{i_1}\wedge \underline{i_2-n}\wedge \underline{i_3}\wedge \dots+ \dots\end{aligned} \end{align} where the usual rules of a wedge product are used to rewrite the right-hand side in terms of the basis of Fock space. Note that since $i_{j+1}=i_j-1$ for all $j\gg 0$, the righthand sum is in fact finite. Then $\alpha_n$ and $\alpha_{-n}$ are adjoint with respect to $\langle\cdot | \cdot \rangle$ and preserve the charge. While $e$ acts by zero on $V$, it must act by the identity on $\Lambda^{\infty/2}V$ for (\ref{action}) to define a representation. We say the representation on Fock space has {\it central charge} $1$. A {\it vertex operator} is an expression of the form  \begin{align*} \psi^+(c):=\textrm{exp}\left(\sum_{n> 0} c_n\frac{\alpha_n}{n} \right)\hspace{5pt}\textrm{    or    }\hspace{5pt}\psi^-(c)=\textrm{exp}\left(\sum_{n< 0} c_n\frac{\alpha_n}{n} \right).\end{align*} Strictly speaking, $\psi^-(c)$ does not act on Fock space, since applying it to a basis vector gives an infinite linear combination. If we complete Fock space, then for appropriate choices of constants $\{c_n\}$ the expressions converge. 

Define an analytic function of $z$ depending on a partition $\lambda$ as follows:  \begin{align}\label{edef}{\bf e}(\lambda,z):=\sum_{i= 1}^{\infty} e^{(\lambda_i-i+\frac{1}{2})z}.\end{align} Observe that since $\lambda_i=0$ eventually, this function admits a meromorphic continuation to all $\C$ via the geometric series formula. Then we may define a function ${\bf p}_k(\lambda)$ on partitions as the $k$th Taylor coefficient in the expansion of ${\bf e}(\lambda,z)$ about $z=0$:  \begin{align*} {\bf e}(\lambda,z)=\frac{1}{z}+\sum_{k\geq 1} {\bf p}_k(\lambda)\frac{z^k}{k!}.\end{align*} The functions ${\bf p}_k(\lambda)$ play an essential role in the Gromov-Witten or Hurwitz theory of an elliptic curve, see Eskin-Okounkov \cite{eo1} and Okounkov-Pandharipande \cite{op}. Following these references, we define an operator which acts diagonally in the $v_\lambda$ basis of $\mathbb{L}$ by \begin{align*}\mathcal{E}_0(z)v_\lambda= {\bf e}(\lambda,z)v_\lambda.\end{align*}
%This action naturally extends to the whole Fock space, acting diagonally in the fermionic basis by an eigenvalue analogous to (\ref{edef}). Then we have the commutation relation \begin{align}\mathcal{E}_0(z)\circ S = e^z(S\circ \mathcal{E}_0(z))\end{align} with the shift operator.
There is an expression for the action of $\mathcal{E}_0$ on $\mathbb{L}$ in terms of vertex operators, see e.g. formula (33) of \cite{eo2}:  \begin{align} \label{caw}\mathcal{E}_0(\ln x)=\frac{1}{x^{1/2}-x^{-1/2}}\,[y^0]\,\textrm{exp}\left(\sum_{n\neq 0} \frac{(xy)^n-y^n}{n}\alpha_n \right).\end{align}

The {\it $N$-quotients} $\lambda^{r/N}$ for $r=0,\dots,N-1$ of the partition $\lambda$ are constructed as follows: Define the subset \begin{align*}T_r(
\lambda):=\{\lambda_i-i+\tfrac{1}{2}\,:\,\lambda_i-i\equiv r\textrm{ mod }N\}.\end{align*} Then wedge the elements of $T_r(\lambda)$ together to produce an element $(v_\lambda)_r\in\textstyle\Lambda^{\infty/2} V_r$ where $V_r$ is the subspace of $V$ defined by \begin{align*}V_r:=\textrm{span}\,\{\underline{i}\,:\,i\in \tfrac{1}{2}+r+N\Z\}.\end{align*} Then there is a unique integer $C_r(\lambda)$, which we call the {\it $r/N$-charge} of $\lambda$ such that the isomorphism $V_r\rightarrow V$ sending \begin{align}\label{isor}\underline{i}\,\mapsto\, \underline{\tfrac{1}{N}(i-r-\tfrac{1}{2})+\tfrac{1}{2}}\end{align} sends $(v_\lambda)_r$ to a basis element of Fock space of charge $C_r(\lambda)$. Shifting by $-C_r(\lambda)$ produces a basis vector $v_{\lambda^{r/N}}$ of $\mathbb{L}$, which defines the $r$th quotient partition $\lambda^{r/N}$. We say $\lambda$ is {\it $N$-decomposable} if and only if the charges $C_r(\lambda)=0$ for all $r$. Equivalently, the Young diagram of $\lambda$ admits a decomposition into $N$-rim hooks, also called $N$-ribbons, see Section 3 of \cite{engel} for more details. More generally, if we remove as many $N$-ribbons as possible from $\lambda$, the remainder is the $N$-core $\lambda^{\textrm{mod }N}$, which is determined uniquely by the charges $C_r(\lambda)$. These integers form a point in the $A_{N-1}$-lattice \begin{align*}A_{N-1}:=\left\{(C_0,\dots,C_{N-1})\in \Z^N\,:\,\textstyle\sum C_i=0\right\}.\end{align*}

We may collect the inverses of the isomorphisms $V_r\rightarrow V$ from (\ref{isor}) into an isomorphism $\Lambda^{\infty/2}V \rightarrow (\Lambda^{\infty/2}V)^{\otimes N}$ which breaks a half-infinite wedge into the sub-wedges lying in a given congruence class mod $N$. Decomposing by charge, we have an isomorphism  \begin{align}\label{tensordecomp}\begin{aligned} \mathbb{L}\otimes \C[\Z] & \rightarrow \mathbb{L}^{\otimes N}\otimes \C[\Z^N] \\ 
v_{\lambda}\otimes C &\mapsto \pm (v_{\lambda^{0/N}}\otimes \cdots \otimes v_{\lambda^{N-1/N}})\otimes (C_0,\dots,C_{N-1})\end{aligned}\end{align} such that $\sum C_r=C$. We will clarify the sign later. We define an operator $\mathcal{E}_0^{r/N}(z)$ which can be thought of as only acting on the $r$th tensor factor in a manner analogous to $\mathcal{E}_0(z)$ \begin{align*}\mathcal{E}_0^{r/N}(z)v_\lambda :={\bf e}(\lambda^{r/N}\!,z)v_\lambda.\end{align*}

\begin{remark} It is possible to express $\mathcal{E}_0^{r/N}(z)$ in terms of $\mathcal{E}_0(z)$ by observing  \begin{align*} \sum_{j\in T_r(\lambda)} e^{jz} &=\frac{1}{N}\sum_{s\,\,(N)} \zeta^{-(r+\frac{1}{2})s}\mathcal{E}_0(z+\tfrac{2\pi i s}{N}) \\ &=e^{\left(NC_r+\tfrac{2r+1-N}{2}\right)z}\mathcal{E}_0^{r/N}(Nz),\end{align*} but it is unnecessary to do so, as we will mostly work with the righthand side of (\ref{tensordecomp}) as opposed to the lefthand side. \end{remark}

%To express this operator in terms of the Heisenberg algebra, we introduce operators $\alpha^{r/N}_k$ which act by $\alpha_k$ on the $r$th tensor factor of the righthand side of (\ref{tensordecomp}) and act trivially on all other factors. Then by formula (\ref{caw}) \begin{align}\mathcal{E}_0^r(\ln x) = \frac{x^{NC_r+\tfrac{2r+1-N}{2}}}{x^{N/2}-x^{-N/2}}\,[y^0]\,\textrm{exp}\left(\sum_{n\neq 0} \frac{y^{-n}-(x^Ny)^{-n}}{n}\alpha^{r/N}_n \right).\end{align}

The operator $q^H$ on Fock space itself decomposes as a tensor product of operators with respect to (\ref{tensordecomp}) as \begin{align*}q^H = (q^{NH_0}\otimes \cdots\otimes q^{NH_{N-1}})\,\otimes\,q^{Q(C_0,\,\dots,\,C_{N-1})}\end{align*} where $Q$ is the quadratic polynomial \begin{align}\label{form}Q(C_0,\,\dots,\,C_{N-1})=\frac{N||\vec{C}||^2}{2}+\sum_{r=0}^{N-1} \left(r+\tfrac{1-N}{2}\right) \cdot C_r.\end{align} When restricted to $\mathbb{L}$, the charge zero subspace, $Q$ gives the formula for the size of the $N$-core of a partition in terms of its charges \cite{cores}. 

Define {\it shift operators $S_r$} which act by $C_r\mapsto C_r+1$ and trivially on all other $r/N$-charges. Then applying a composition of shift operators and their inverses corresponds to changing the $N$-core of $\lambda$ without changing its quotients. In particular, define \begin{align*}S_\eta:=\!\!\prod_{r\textrm{ mod }N}\!\! S_r^{C_r(\eta)}\end{align*} where $\eta$ is some $N$-core. Thus if $\lambda^{\textrm{mod }N}=\eta$, the partition associated to $S_\eta^{-1}v_\lambda$ is $N$-decomposable. Note that $S_\eta$ and $S_\eta^{-1}$ are adjoint isometries with respect to $\langle \cdot \big{|}\cdot\rangle$. Let $\mathfrak{R}_N$ be the {\it rotation operator} which cyclically rotates the factors of the tensor decomposition (\ref{tensordecomp}), acting both on $\mathbb{L}^{\otimes N}$ and the charge lattice $\Z^N$. Then, $\mathfrak{R}_N$ cyclically permutes the $r/N$-charges, so that \begin{align*}C_r(\mathfrak{R}_N(\eta))=C_{r-1}(\eta).\end{align*}

We now recall some definitions from \cite{engel}:

\begin{definition} Let $\eta$ be an $N$-core. The {\it $\eta$-weight} of a partition is \begin{align*}{\bf w}_{N,\eta}(\lambda):= \twopartdefotherwise{\pm u_\eta\displaystyle\frac{\langle 1\rangle \langle N-1\rangle}{\langle 0 \rangle^2}}{\lambda^{{\rm mod}\,N}= \eta}{0}\end{align*} Here $\langle a\rangle$ denotes the product of the hook lengths of $\lambda$ congruent to $a$ mod $N$, and $u_\eta$ is the unique constant making ${\bf w}_{N,\eta}(\eta)=1$. The sign is defined so that Theorem \ref{big} below is true. \end{definition}

In \cite{engel}, we defined an enlargement of the algebra of shifted symmetric functions to be the ring functions on partitions generated by the Taylor coefficients of $\mathcal{E}_0(z+\tfrac{2\pi ri}{N})$. Here, we give a slightly different definition:

\begin{definition}\label{cyclotomic} The {\it cyclotomic enlargement of the algebra of shifted symmetric functions} is the ring of functions on partitions $$\Lambda_N:=\C[{\bf p}_k(\lambda^{r/N})].$$ \end{definition} 

\begin{remark} This definition differs from that in \cite{engel}, but not when one restricts to the set of all $\lambda$ with the same $N$-core---Definition \ref{cyclotomic} gives a ring of functions which is insensitive to applying the shift operators $S_r$, unlike the functions ${\bf p}_k^r(\lambda)$ from \cite{engel}. \end{remark} 

Finally, we define the following generalization of the Bloch-Okounkov bracket \cite{bo} ($N=1$ and $\eta=\emptyset$) and the pillowcase bracket ($N=2$ and $\eta=\emptyset$) of Eskin-Okounkov \cite{eo2}.

\begin{definition} Let $\eta$ be an $N$-core. Let ${\bf f}$ be a function from the set of partitions to $\C$. Its {\it $\eta$-bracket} is \begin{align*}\langle {\bf f}\rangle_{{\bf w}_{N,\eta}}:=\frac{\sum_\lambda {\bf f}(\lambda){\bf w}_{N,\eta}(\lambda)q^{|\lambda|/N}}{\sum_\lambda {\bf w}_{N,\emptyset}(\lambda)q^{|\lambda|/N}}.\end{align*} \end{definition}

This definition also differs slightly from \cite{engel}, where the exponent of $q$ is not divided by $N$. We may now state the main theorem of this paper:

\begin{theorem}\label{modularity} The image of $\langle \,\cdot\, \rangle_{{\bf w}_{N,\eta}}\,:\,\Lambda_N\rightarrow \C[[q^{1/N}]]$ lies in the ring of quasimodular forms for $\Gamma(N)$.
%Furthermore, ${\rm wt}\,\langle {\bf f}\rangle_{{\bf w}_{N,\eta}}={\rm wt}({\bf f})+|\eta|$.
\end{theorem}

The relevance of this theorem lies in that fact that when $N=1,2,3,4,6$ the brackets ${\rm wt}\,\langle {\bf f}\rangle_{{\bf w}_{N,\eta}}$ fully determine the Hurwitz theory of the elliptic orbifold $E/\langle \zeta_N\rangle$. For instance, as shown in Section 2 of \cite{engel} this theorem implies that generating functions of triangulations with fixed curvatures are quasimodular forms. 

\begin{definition} The {\it $\eta$-weighted $n$-point function} is \begin{align}\label{npoint}F_{N,\eta}\begin{pmatrix} z_1 & \dots & z_n \\ r_1 & \dots & r_n\end{pmatrix}:=\left\langle {\bf e}(z_1,\lambda^{r_1/N})\dots {\bf e}(z_n,\lambda^{r_n/N})\right\rangle_{{\bf w}_{N,\eta}}.\end{align} \end{definition}

\begin{remark} To prove Theorem \ref{modularity}, it suffices to show that the Taylor coefficients $[z_1^{k_1}\dots z_n^{k_n}]$ of this function are quasimodular forms for $\Gamma(N)$---this in turn implies that the $\eta$-bracket of any monomial in the ${\bf p}_k(\lambda^{r/N})$'s is quasimodular, and these form a $\C$-basis for $\Lambda_N$. The strategy (now common for such questions) is to explicitly compute (\ref{npoint}) by realizing it as the trace of a product of operators on Fock space. We are nearly already able to do so, as the $\mathcal{E}_0^{r/N}(z)$ operators act diagonally in the fermionic basis of $\mathbb{L}$, as does $q^{H/N}$. But we do not yet possess a vertex operator whose diagonal entry for $v_\lambda$ is ${\bf w}_{N,\eta}(\lambda)$. This main technical obstacle is resolved in Theorem \ref{big}. \end{remark}

Define an enlargement of the Heisenberg algebra action as follows. Let $\alpha_n^{r/N}$ be the operator on Fock space induced by the following action on $V$: \begin{align*}\alpha_n^{r/N}(\underline{i})=\twopartdefotherwise{\underline{i-Nn}}{i-\tfrac{1}{2}\equiv r\textrm{ mod }N}{0}\end{align*} 
%On Fock space, these operators satisfy the commutation relations  \begin{align} \alpha_n^{r/N}\alpha_m &= \alpha_m\alpha_n^{r+m/N}\hspace{10pt}\textrm{if }(m,N)=1 \\ [\alpha_n^{r/N},\alpha_m^{s/N}]&=n\delta_{n,-m}\delta_{r,s}.\end{align} 
We now have the machinery to describe the {\it $N$-core operators} $\mathfrak{W}_{N,\eta}$ whose diagonal entries are the $\eta$-weights:

\begin{theorem}\label{big} Let $\eta$ be an $N$-core with charges $C_r:=C_r(\eta)$. Define  \begin{align*} \psi_\eta:=&{\rm exp}\,\bigg(\!\!\sum_{\substack{n>0 \\ r\textrm{ mod }N}}\!\!\! (C_{r-1}-C_r)\frac{\alpha_n^{r/N}}{n} \bigg), \\ \mathfrak{W}_N:=&{\rm exp}\left( \sum_{n\in \Z\backslash N\Z} \frac{\alpha_n}{n}\right). \end{align*}  Define the operator \begin{align*}\mathfrak{W}_{N,\eta}:=S_\eta\circ \psi_{\eta}^{-T} \circ \mathfrak{W}_N \circ \psi_{\mathfrak{R}_N^{-1}(\eta)} \circ S_\eta^{-1}.\end{align*} Then we have \begin{align*}\langle v_\lambda \,\big{|}\,\mathfrak{W}_{N,\eta}\, \big{|}\, v_\lambda\rangle={\bf w}_{N,\eta}(\lambda):= \twopartdefotherwise{\pm u_\eta\displaystyle\frac{\langle 1\rangle \langle N-1\rangle}{\langle 0 \rangle^2}}{\lambda^{{\rm mod}\,N}= \eta}{0}\end{align*} Here $\langle a\rangle$ denotes the product of the hook lengths of $\lambda$ congruent to $a$ mod $N$, and $u_\eta$ is the unique constant making ${\bf w}_{N,\eta}(\eta)=1$. \end{theorem} 

\begin{proof} For the remainder of the proof, let $v_\mu:=S_\eta^{-1}v_\lambda$. The theorem is equivalent to showing that ${\bf w}_{N,\eta}(\lambda)=\langle v_\mu \,|\, \mathfrak{O}_{N,\eta}\,|\,v_\mu\rangle$ where \begin{align*}\mathfrak{O}_{N,\eta}:={\rm exp}\bigg(\!\!\sum_{\substack{n<0 \\ r\textrm{ mod }N}}\!\!\! (C_{r-1}-C_{r})\frac{\alpha_n^{r/N}}{n} \bigg)\circ \mathfrak{W}_N\circ{\rm exp}\bigg(\!\!\sum_{\substack{n>0 \\ r\textrm{ mod }N}}\!\!\! (C_{r}-C_{r+1})\frac{\alpha_n^{r/N}}{n} \bigg).\end{align*} For any given $\mu$, we may compute a matrix entry of $\mathfrak{O}_{N,\eta}$ by taking the determinant of a minor of the matrix entries of $\mathfrak{O}_{N,\eta}$ acting on a truncation \begin{align*}V^{N\ell}=\textrm{span}\{\underline{i}\,:\,|i|\leq N\ell-\tfrac{1}{2}\}.\end{align*} This step is justified because $\mathfrak{O}_{N,\eta}$ is a product of upper and lower unitriangular matrices. It is useful to re-index the basis vectors by setting \begin{align*}e_i = \underline{i+\tfrac{1}{2}-N\ell}\end{align*} so that $i\in\{0,\dots,N\ell-1\}$. On $V^{N\ell}$ the action of $\alpha_n$ is \begin{align*}\alpha_n(e_i)=\alpha_1^n(e_i)=	\left\{
		\begin{array}{ll}
			\underline{i-n} & \mbox{if } i-n\geq 0 \\
			0 & \mbox{otherwise }
		\end{array}
	\right.\end{align*} and the action of $\alpha^{r/N}_n$ is given by \begin{align*}\alpha^{r/N}_n(e_i)=\alpha_1^{r/N}(e_i)=\twopartdefotherwise{e_{i-Nn}}{i-\tfrac{1}{2}\equiv r\textrm{ mod }N\textrm{ and }i-Nn\geq 0}{0}\end{align*} We now compute the matrix entry $\langle e_i\,|\,\mathfrak{O}_{N,\eta}\,|\,e_j\rangle$. We have  \begin{align*} {\rm exp}\bigg(\!\!\sum_{\substack{n>0 \\ r\textrm{ mod }N}}\!\!\! (C_{r}-C_{r+1})\frac{\alpha_n^{r/N}}{n} \bigg) e_j &= {\rm exp}\bigg(\sum_{n>0} (C_{j}-C_{j+1})\frac{\alpha_n^{j/N}}{n} \bigg)e_j \\ &=\sum_{k\geq 0} e_{j-Nk}\,[y^{Nk}](1-y^N)^{C_{j+1}-C_{j}} \end{align*} because ${\rm exp}\left(\sum_{n>0}\frac{y^{Nn}}{n}\right)=(1-y^N)^{-1}$. We may similarly compute the action of the leftmost term in $\mathfrak{O}_{N,\eta}$ on $e_i$ to conclude that \begin{align*}\langle e_i\,|\,\mathfrak{O}_{N,\eta}\,|\,e_j\rangle=[x^iy^j]\,\frac{(1-y^N)^{C_{j+1}-C_{j}}} {(1-x^N)^{C_{i}-C_{i-1}}}\sum_{k,l\geq 0} x^ky^l  \left\langle e_k\,|\, \mathfrak{W}_N\,|\,e_l\right\rangle.\end{align*} The matrix entries of $\mathfrak{W}_N$ on $V^{N\ell}$ were computed in \cite{engel} to be: \begin{align*}\sum_{k,l\geq 0} x^ky^l\left\langle e_k\,|\, \mathfrak{W}_N\,|\,e_l\right\rangle=\frac{1}{1-xy}\cdot\frac{1-x}{(1-x^N)^{1/N}}\cdot \frac{(1-y^N)^{1/N}}{1-y}.\end{align*}
	
Expanding by the binomial theorem, we conclude that $\langle e_i\,|\,\mathfrak{O}_{N,\eta}\,|\,e_j\rangle$ is the $[x^iy^j]$ coefficient of  \begin{align*} 
%[x^iy^j]\,\frac{1}{1-xy}\cdot\frac{1-x}{1-y} \cdot (1-x^N)^{-C_s(\eta)-1/N}\cdot (1-y^N)^{-C_t(\eta)+1/N} \\ 
\sum_{a,b,c\geq 0}(-1)^{a+b}{C_{i-1}-C_{i}-\tfrac{1}{N}  \choose a}{C_{j+1}-C_{j}+\tfrac{1}{N}  \choose b} (x^{aN}-x^{aN+c+1})y^{bN+c}. \end{align*} Observe that $\langle e_i\,|\,\mathfrak{O}_{N,\eta}\,|\,e_j\rangle=0$ unless either $i\equiv 0\,(N)$ or $i\equiv j+1\,(N)$. Let  \begin{align*}{\bf b}_\eta(i):=&\prod_{r = 1}^i\frac{\{r+N(C_i-C_{i-1})\,\big{|}\,r\equiv 1\textrm{ mod }N\}}{\{r\,\big{|}\,r\equiv 0\textrm{ mod }N\}}\hspace{10pt}\textrm{and} \\ {\bf c}_\eta(j):=&\prod_{r = 1}^j   \frac{\{r+N(C_j-C_{j+1})\,\big{|}\,r\equiv -1\textrm{ mod }N\}}{\{r\,\big{|}\,r\equiv 0\textrm{ mod }N\}}.\end{align*} Simple manipulations with binomial coefficients show that \begin{align*}\langle e_i\,|\,\mathfrak{O}_{N,\eta}\,|\,e_j\rangle={\bf b}_\eta(i){\bf c}_\eta(j)\end{align*} whenever $i\equiv 0\,(N)$ and $j\not\equiv -1\,(N)$. We now compute the ``interesting" matrix entries for which $i\equiv j+1\,(N)$. Let $a= \lfloor \frac{i}{N}\rfloor$ and  $b= \lfloor \frac{j}{N}\rfloor$ and assume for convenience that $i\not\equiv 0$ and $b>a$. By applying the identity   \begin{align*} {X \choose 0}{-X \choose b-a}+\cdots +{X \choose a}{-X  \choose b} ={b+X \choose a}{a-X \choose b} \end{align*} to $X = C_{i-1}-C_{i}-\frac{1}{N}$, we compute \begin{align*}\langle e_i\,|\, \mathfrak{O}_{N,\eta} \,|\, e_j \rangle = (-1)^{a+b+1}{b+X \choose a}{a-X \choose b}  = \frac{{\bf b}_\eta(i){\bf c}_\eta(j)}{NC_{[j]}+j-i-NC_{[i]}}.\end{align*} Similar computations give the same answer when $i\equiv 0\,(N)$, $j\equiv -1\,(N)$. Thus, we have computed all the entries:  \begin{align}\label{entries} \frac{\langle e_i\,|\, \mathfrak{O}_{N,\eta} \,|\, e_j \rangle}{{\bf b}_\eta(i){\bf c}_\eta(j)}=\left\{
	\begin{array}{lll}
		(NC_{[j]}+j-i-NC_{[i]})^{-1} & \mbox{if } i\equiv j+1\,(N) \\
		1 & \mbox{if } i\equiv 0,\, j\not\equiv -1\,(N) \\
		0 & \mbox{otherwise.}
	\end{array}
\right.\end{align}

We have that \begin{align*}\langle v_\lambda \,|\,\mathfrak{W}_{N,\eta}\,|\,v_\lambda\rangle= \langle v_\mu \,|\,\mathfrak{O}_{N,\eta}\,|\,v_\mu\rangle= \det \,\langle e_i \,|\, \mathfrak{O}_{N,\eta} \,|\, e_j\rangle_{i,\,j\in\{\mu_k-k+N\ell\}}.\end{align*} This minor of $\langle e_i\,|\,\mathfrak{O}_{N,\eta}\,|\,e_j\rangle$ can be put into $N\times N$ block form \begin{align}\label{block}\begin{pmatrix} * &  *  & \cdots & * & * \\ * & 0 & \cdots & 0 & 0 \\ 0 & * & \cdots & 0 & 0 \\ \vdots  & \vdots & & \vdots & \vdots \\ 0  & 0 & \cdots & * & 0\end{pmatrix}\end{align} by collecting the indices which lie in each congruence class mod $N$. More precisely, the $(r,s)$-block is the matrix \begin{align*}\langle e_{\mu_i-i+\ell N} \,|\, \mathfrak{O}_{N,\eta}\,|\,e_{\mu_j-j+\ell N}\rangle_{\mu_i-i\equiv r}^{\mu_j-j\equiv s}.\end{align*} This requires reordering the rows and columns by the same permutation, and does not change the determinant. Most blocks are zero according to (\ref{entries}), whereas the $(0,s)$-blocks for $s\neq N-1$ are all rank $1$ matrices---up to rescaling rows of the above matrix by ${\bf b}_\eta(\mu_i-i+\ell N)$ and the columns by ${\bf c}_\eta(\mu_j-j+\ell N)$, these $(0,s)$-blocks will have all entries equal to $1$ by (\ref{entries}).

The width of the $r$th row of blocks of (\ref{block}) is \begin{align*}n_r:=\#\{\mu_i-i+\ell N\geq 0\,:\,\mu_i-i\equiv r\textrm{ mod }N\}.\end{align*} We now apply an argument from \cite{engel} to show all $n_r$ are equal. We first remark that for $0\leq r\leq N-2$, the determinant of (\ref{block}) vanishes unless \begin{align*}n_r-1 \leq n_{r+1} \leq n_r.\end{align*} The upper bound is immediate because the $(r+1,r)$-block is the only non-zero entry in the $(r+1)$th row. On the other hand, the only other non-zero entry in the $r$th column is the $(0,r)$-block, which has rank one. This gives the lower bound. Similar logic applied to the $(0,N-1)$-block implies that \begin{align*}n_{N-1}\leq n_0\leq n_{N-1}+1.\end{align*} Next, observe there is at most one $r\in\{0,\dots,N-2\}$ such that $n_r>n_{r+1}$ as otherwise the $(0,N-1)$ block cannot have the correct proportions. Therefore \begin{align*}Nn_0\geq \sum n_r \geq Nn_0-(N-1).\end{align*} Since $n_r = C_r(\mu)+\ell$, we have $\sum n_r = N\ell$ is divisible by $N$. Therefore $\langle v_\lambda \,|\,\mathfrak{W}_{N,\eta}\,|\,v_\lambda\rangle=0$ unless all blocks in (\ref{block}) are square , i.e. $C_r(\mu)=0$ for all $r$. Equivalently, $\lambda$ must have $N$-core equal to $\eta$. This verifies the theorem in the case where ${\bf w}_{N,\eta}(\lambda)=0$.

If $\lambda^{\textrm{mod }N}=\eta$, the determinant is, up to sign, the product of the determinants of the $(r+1,r)$-blocks. We apply the Cauchy determinant formula \begin{align*}\det\left(\frac{1}{x_k-y_\ell}\right) = \frac{\prod_{k<l}(x_k-x_\ell)(y_\ell-y_k) }{\prod_{k,l} (x_k-y_\ell)}\end{align*} to each block to conclude  \begin{align}\label{superprod}\begin{aligned} &\langle v_\lambda\,|\,\mathfrak{W}_{N,\eta}\,|\,v_\lambda\rangle = \pm\prod_{i,j} \,{\bf b}_\eta(\mu_i-i+N\ell) {\bf c}_\eta(\mu_j-j+N\ell)\cdot \\ &\frac{\{\mu_j-\mu_i+i-j\equiv 0\textrm{ mod }N\}^2}{\{\mu_j-\mu_i+i-j+NC_{\mu_j-j}-NC_{\mu_i-i}\equiv \pm 1\textrm{ mod }N\}}.\end{aligned}\end{align} Roughly, for each $r\equiv \pm 1\textrm{ mod }N$ the number of times $r$ appears in the numerator will be the number of $r$-hooks of $\lambda$, whereas for each $r\equiv 0 \textrm{ mod }N$ the number of times $r$ appears in the denominator will be twice the number of $r$-hooks. But this heuristic is not accurate for small $r$. We show this ``error" exactly accounts for the factor $u_\eta$ which appears in ${\bf w}_{N,\eta}(\lambda).$ This requires an interesting combinatorial argument which we relegate to Appendix A. \end{proof}

 We now return to the analysis of the $\eta$-weighted $n$-point function defined in (\ref{npoint}). It is convenient to change variables by setting $x_i=e^{z_i}$. Then by Theorem \ref{big} we have  \begin{align*} F_{N,\eta}\begin{pmatrix} z_1 & \dots & z_n \\ r_1 & \dots & r_n\end{pmatrix}=\frac{tr_{\mathbb{L}}\,\,q^{H/N}\mathcal{E}_0^{r_1/N}(\ln x_1)\dots \mathcal{E}_0^{r_n/N}(\ln x_n)\mathfrak{W}_{N,\eta}}{tr_{\mathbb{L}}\,\,q^H\mathfrak{W}_N}.\end{align*} For notational compactness, we write the above function as $F_{N,\eta}(z,r)$ where $z=(z_1,\dots,z_n)$ and $r=(r_1,\dots,r_n)$. Since all operators but $\mathfrak{W}_{N,\eta}$ act diagonally in the fermionic basis, the only contribution to the trace in the numerator is from the subspace \begin{align*}\mathbb{L}^{\otimes N}\otimes (C_r(\eta))\end{align*} which is spanned by the $v_\lambda$ for which $\lambda^{\textrm{mod }N}=\eta$. First, we move the left-hand shift operator in \begin{align*}\mathfrak{W}_{N,\eta}=S_\eta\circ \mathfrak{O}_{N,\eta}\circ S_{\eta}^{-1}\end{align*} further left through the remaining operators, by the commutation relations  \begin{align*} \mathcal{E}_0^{r/N}(\ln x)\circ S_\eta &= S_\eta\circ \mathcal{E}_0^{r/N}(\ln x) \\ \langle v_\lambda \,|\, q^{H/N}\circ S_\eta
%\big{|}_{\mathbb{L}^{\otimes N}\otimes \{0\}}
&=
%\frac{q^{Q(C_0(\eta)+C_0,\,\dots,\,C_{N-1}(\eta)+C_{N-1})}}{q^{Q(C_0,\,\dots,\,C_{N-1})}}
q^{|\eta|/N}\langle v_\lambda\,|\,S_\eta \circ q^{H/N}\end{align*} assuming that $\lambda^{\textrm{mod }N}=\eta$.
Having moved the shift operators to the extremities, we may remove them without changing the trace:
 \begin{align}\label{traceform}\begin{aligned} F_{N,\eta}(z,r) =q^{|\eta|/N}\frac{tr_{\mathbb{L}}\,\,q^{H/N}\mathcal{E}_0^{r_1/N}(\ln x_1)\dots \mathcal{E}_0^{r_n/N}(\ln x_n)\mathfrak{O}_{N,\eta}}{tr_{\mathbb{L}}\,\,q^H\mathfrak{W}_N}. \end{aligned}\end{align}

The only elements $v_\lambda$ which contribute to this trace are those for which $\lambda$ is $N$-decomposable. Under the isomorphism (\ref{tensordecomp}), they correspond to the fermionic basis of the subspace with all $r/N$-charges equal to zero: \begin{align*}\mathbb{L}^{\otimes N}=\mathbb{L}^{\otimes N}\otimes (0,\dots,0).\end{align*}

\begin{remark} We finally clarify the ambiguous sign in (\ref{tensordecomp}). Essentially, we would like the operators $\alpha_n^{r/N}$ to act solely on the $r$th tensor factor on the right side of (\ref{tensordecomp}) so that we may express \begin{align*} \mathcal{E}_0^{r/N}(z)=1\otimes \cdots \otimes \mathcal{E}_0(z)\otimes \cdots\otimes 1\end{align*} using the formula (\ref{caw}), but with $\alpha_n^{r/N}$ replacing $\alpha_n$. But the action of \begin{align}\label{otheralpha} 1\otimes \cdots \otimes \alpha_n \otimes \cdots\otimes 1\end{align} and $\alpha_n^{r/N}$ do not agree if we always put a positive sign in (\ref{tensordecomp})---applying $\alpha_n^{r/N}$ to an element $v_\lambda$ gives a signed linear combination of partitions, whose signs do in fact depend on more data than the quotient partition $\lambda^{r/N}$ alone. But we may choose a compatible collection of signs in the isomorphism (\ref{tensordecomp}) such that the action of $\alpha_n^{r/N}$ on $\mathbb{L}$ agrees with the action of (\ref{otheralpha}) on $\mathbb{L}^{\otimes N}$. Such a collection of signs can be constructed by taking a sequence of truncations $V^{[\ell]}$ and applying the more pedestrian isomorphism $$\Lambda^kV^{[\ell]}\rightarrow\bigoplus_{\sum i_r=k} \Lambda^{i_r} V_r^{[\ell]}$$ in a compatible manner as $\ell\rightarrow \infty$. Then the two actions agree. \end{remark}

The operators $\alpha_n^{r/N}$ on $\mathbb{L}^{\otimes N}$ satisfy the commutation relations \begin{align*}[\alpha_n^{r/N},\alpha_m^{s/N}]= n\delta_{n,-m}\delta_{rs}.\end{align*} Having corrected signs, all of the operators in (\ref{traceform}) are expressible terms of the $\alpha_n^{r/N}$ with the lone exception of $\mathfrak{W}_N$ appearing in $\mathfrak{O}_{N,\eta}$. Thus, it would be advantageous to find an operator $\overline{\mathfrak{W}}_N$ acting on $\mathbb{L}^{\otimes N}$ such that  \begin{align*}\langle v_\lambda \,\big{|}\, \mathfrak{W}_N\,\big{|}\,v_\lambda \rangle= \langle v_\lambda \,\big{|}\, \overline{\mathfrak{W}}_N\,\big{|}\,v_\lambda \rangle \end{align*} for any $N$-decomposable $\lambda$. From the proof of Theorem \ref{big}, the only relevant matrix entries of $\langle e_i\,\big{|}\,\mathfrak{W}_N\,\big{|}\,e_j\rangle$ in the computation of the the $\eta$-weight are those for which $i\equiv j+1\textrm{ mod }N$. Thus, we make the educated guess that \begin{align*}\overline{\mathfrak{W}}_N=(\mathfrak{v}_0\otimes\cdots \otimes \mathfrak{v}_{N-1}) \circ \mathfrak{R}_N\end{align*} with $\mathfrak{v}_i$ acting on the $i$th tensor factor of $\mathbb{L}^{\otimes {N}}$. While the action of $\mathfrak{R}_N$ simply rotates the factors of $\mathbb{L}^{\otimes N}$, the manner in which it acts on the $N$-decomposable subspace of $\mathbb{L}$ is more subtle---it sends  \begin{align}\label{issue} \mathfrak{R}_N\,:\,\underline{i}\mapsto \twopartdef{\underline{i+1}}{i\not\equiv -1\textrm{ mod }N}{\underline{i+1-N}}{i\equiv -1 \textrm{ mod }N.} \end{align}

Assume $r\not\equiv 0\textrm{ mod }N$. Consider the matrix entries (\ref{entries}) of $\mathfrak{W}_N$ lying in the $(r,r-1)$-block, i.e. those for which $i \equiv r$ and $j\equiv r-1$. By composing with the isomorphisms $V_r\mapsto V$ and $V_{r-1}\mapsto V$, then applying $\mathfrak{R}_N$ to $\underline{j}$, we can determine the desired matrix entries of the hypothetical operator $\mathfrak{v}_r$: \begin{align*}\langle e_a \,\big{|}\, \mathfrak{v}_r \,\big{|}\,e_b\rangle=\frac{{\bf b}(Na+1){\bf c}(Nb)}{N(b-a)-1}=\frac{(-1)^{a+b}}{N(b-a)-1}{-1-\tfrac{1}{N}\choose a}{-1+\tfrac{1}{N} \choose b}.\end{align*} Here $a=\lfloor i/N\rfloor$ and $b=\lfloor j/N\rfloor$ play the same role that they did in the proof of the Theorem. When $r\equiv 0\textrm{ mod }N$, we have  the unusual case in (\ref{issue}), but the values of ${\bf b}(i)$ and ${\bf c}(j)$ are also unusual. We would like $\mathfrak{v}_0$ to have matrix entries \begin{align*}\langle e_a \,\big{|}\,\mathfrak{v}_0 \,\big{|}\,e_b\rangle =\frac{{\bf b}(Na){\bf c}(N(b+1)-1)}{N(b-a+1)-1}.\end{align*}

%= \frac{(-1)^{a+b+1}N(b+1)}{N(b+1-a)-1}{-\tfrac{1}{N}\choose a}{-1+\tfrac{1}{N} \choose b+1}.\end{align}

Rephrasing, when $r\not\equiv 0$, we would like to find an operator $\mathfrak{v}_r$ such that \begin{align*}G_r(x,y) := \sum_{a,b,\geq 0} \langle e_a\,\big{|}\,\mathfrak{v}_r\,\big{|}\,e_b\rangle \,x^ay^b\end{align*} satisfies the differential equation \begin{align*}(Ny\partial_y-Nx\partial_x-1)G_r(x,y)=-(1-x)^{-1-\frac{1}{N}}(1-y)^{-1+\frac{1}{N}}.\end{align*} There is a similar type of equation when $r\equiv 0$. Rather incredibly, these equations have exact closed form solutions:  \begin{align*} G_r(x,y) = \frac{1}{1-xy}\cdot \frac{(1-y)^{1/N}}{(1-x)^{1/N}} && \textrm{ if }r\not\equiv 0 \\ G_r(x,y) = \frac{1}{1-xy}\cdot \frac{(1-y)^{1/N-1}}{(1-x)^{1/N-1}} && \textrm{ if }r\equiv 0. \end{align*}

This allows us to reverse engineer the operator $\overline{\mathfrak{W}}_N$:

\begin{proposition}\label{rotator} Define operators on $\mathbb{L}^{\otimes N}$  \begin{align*} \mathfrak{v}_r:=&{\rm exp}\bigg(\sum_{n\neq 0} (\delta_{r0}-\tfrac{1}{N})\frac{\alpha_n^{r/N}}{n}\bigg) \\ \overline{\mathfrak{W}}_N:=&(\mathfrak{v}_0\otimes \cdots\otimes \mathfrak{v}_{N-1})\circ \mathfrak{R}_N. \end{align*} Then, for any $N$-decomposable partition $\lambda$, we have  \begin{align*}\langle v_\lambda\,\big{|}\,\mathfrak{W}_N\,\big{|}\,v_\lambda\rangle =\langle v_{\lambda^{0/N}}\otimes \cdots\otimes v_{\lambda^{N-1/N}} \,\big{|}\,\overline{\mathfrak{W}}_N\,\big{|}\,v_{\lambda^{0/N}}\otimes \cdots\otimes v_{\lambda^{N-1/N}}\rangle.\end{align*} \end{proposition} 

\begin{remark}\label{rotator2} Proposition \ref{rotator} sheds some light on the rather unusual definition of $\mathfrak{W}_{N,\eta}$---given an orthogonal operator such as $\mathfrak{W}_N$ on $\mathbb{L}$, one would normally expect to modify it by \begin{align*}\mathfrak{W}_N\rightarrow A^{-T}\mathfrak{W}_N A\end{align*} but $\mathfrak{W}_{N,\eta}$ is not of this form. Rather, one side has been cyclically rotated. But, replacing $\mathfrak{W}_N$ with $\overline{\mathfrak{W}}_N$, we may pull the operator $\mathfrak{R}_N$ through the righthand side, and the result is  \begin{align*}\overline{\mathfrak{O}}_{N,\eta}:\!&= \psi_\eta^{-T} \circ \overline{\mathfrak{W}}_N \circ \psi_{\mathfrak{R}_N^{-1}(\eta)} \\ &= {\rm exp}\bigg(\!\!\sum_{\substack{n\neq 0 \\ r\textrm{ mod }N}}\!\!\! (C_{r-1}(\eta)-C_r(\eta)+\delta_{r0}-\tfrac{1}{N})\frac{\alpha_n^{r/N}}{n} \bigg) \circ \mathfrak{R}_N \\ &\hspace{130pt}=:\Psi_\eta \circ \mathfrak{R}_N.\end{align*}\end{remark}

By Proposition \ref{rotator} and Remark \ref{rotator2}, we have re-expressed (\ref{traceform}) as the trace over the $N$-decomposable subspace: \begin{align}\label{traceform2}q^{|\eta|/N}\frac{tr_{\mathbb{L^{\otimes N}}}\,\,q^{H/N}\mathcal{E}_0^{r_1/N}(\ln x_1)\dots \mathcal{E}_0^{r_n/N}(\ln x_n)\Psi_\eta\mathfrak{R}_N}{tr_{\mathbb{L^{\otimes N}}}\,\,q^H\mathfrak{W}_N}. \end{align} This is quite attractive, as it {\it nearly} decomposes as the trace of a tensor product of operators, the only exception being the rotation operator $\mathfrak{R}_N$. Regardless, all operators, including the rotation, act nicely in the {\it bosonic basis} of  $\mathbb{L}^{\otimes N}$, also indexed by $N$-tuples of partitions, 
%$|\mu^0\rangle \otimes \cdots \otimes |\mu^{N-1}\rangle$.
% := \prod \alpha_{-\mu^r_i} ^{r/N} \,\,v_\emptyset \otimes \cdots \otimes v_\emptyset.\end{align}
\begin{align*}\mathbb{L}^{\otimes N}=\bigotimes_{r\,(N)}\,\bigotimes_{m=1}^{\infty}\,\,\bigoplus_{n=0}^{\infty}\left(\alpha_{-m}^{r/N}\right)^n v_\emptyset \otimes \cdots \otimes v_\emptyset.\end{align*} Just as the boson-fermion correspondence for $\mathbb{L}$ was the key in computing the $n$-point functions in \cite{eo2,engel}, the same is true here, for $\mathbb{L}^{\otimes N}$. Because of the rotation, the smallest tensor factors on which the operators in (\ref{traceform2}) decompose are \begin{align*}(\mathbb{L}^{\otimes N})_m:=\bigotimes_{r\,(N)}\,\bigoplus_{n=0}^{\infty} \left(\alpha_{-m}^{r/N}\right)^nv_\emptyset \otimes \cdots \otimes v_\emptyset.\end{align*}

Thus, we are led to consider:

\begin{lemma}\label{circularlemma} We have that  \begin{align*} tr_{(\mathbb{L}^{\otimes N})_m}\,\, q^{H/N}{\rm exp}\bigg(\sum_r A_r\frac{\alpha_{-m}^{r/N}}{m}\bigg){\rm exp}\bigg(\sum_r B_r\frac{\alpha_m^{r/N}}{m}\bigg) \mathfrak{R}_N\end{align*} is equal to  \begin{align*}\frac{1}{1-q^{Nm}}\, {\rm exp}\bigg(\sum_{s=1}^N \frac{q^{sm}}{m(1-q^{Nm})}\sum_{r \,(N)} A_rB_{r+s}  \bigg). \end{align*} \end{lemma}

\begin{proof} First we rescale the bosonic basis to produce an orthonormal basis of $(\mathbb{L}^{\otimes N})_m$ so that the trace may be computed as a sum of inner products: \begin{align*}|{\bf n}\rangle=|n_0,\dots,n_{N-1}\rangle:=\frac{1}{\sqrt{m^{|{\bf n}|}{\bf n}!}} \prod_r \left(\alpha^{r/N}_m\right)^{n_r} \!v_\emptyset\otimes \cdots\otimes v_\emptyset\end{align*} where ${\bf n}!=\prod n_r!$ and $|{\bf n}|=\sum n_r$. By the commutation relation, we directly compute that the first quantity in the statement of the lemma equals \begin{align*}\sum_{\bf n} \frac{q^{m|{\bf n}|}}{m^{|{\bf n}|}{\bf n}!}\prod_{r\,(N)} \,\sum_{k \geq 0} A_r^{n_r-k}B_r^{n_{r-1}-k}{n_r \choose k}{n_{r-1} \choose k}k!.\end{align*} We leave it as an exercise for the combinatorially-minded reader to verify that this equals the second quantity in the statement of the lemma. (This equality is non-trivial and was found using Mathematica.) \end{proof}

To apply the lemma to the case at hand, we must normal-order the operator acting on $(\mathbb{L}^{\otimes N})_m$, that is, move all the operators $\alpha_m^{r/N}$ with $m>0$ to the right of (\ref{traceform2}). Again by the commutation relation $[\alpha_m^{r/N},\alpha_{-m}^{s/N}]=m\delta_{rs}$, \begin{align*}{\rm exp}(B \alpha^{r/N}_m) \,{\rm exp}(A \alpha^{s/N}_{-m})={\rm exp}(A \alpha^{s/N}_{-m}) \,{\rm exp}(B \alpha^{r/N}_m)\,{\rm exp}(ABm\delta_{rs}).\end{align*} Normal-ordering and taking the product of the traces from Lemma \ref{circularlemma} over all $m$ gives the trace over $\mathbb{L}^{\otimes N}$. See Appendix B for the resulting very large formula, and some of the steps made to simplify it. Given a congruence class $r$ mod $N$, define the shifted theta function \begin{align*}\vartheta_r(x,q):=(x^{1/2}-x^{-1/2})^{\delta_{r0}} \prod_{m\geq 1} \frac{(1-xq^{N(m-1)+[r]}) (1-x^{-1}q^{N(m-1)+[-r]})}{(1-q^{N(m-1)+[r]}) (1-q^{N(m-1)+[-r]})}\end{align*} where $[r]$ denotes the unique representative of $r$ mod $N$ such that $0<r\leq N$. Combining the terms at the end of Appendix B, we have:

\begin{theorem}\label{theanswer} Let $\eta$ be an $N$-core and define $M(r):=C_{r-1}-C_r+\delta_{r0}-\tfrac{1}{N}.$ In the annular domain \begin{align}\label{domain} q|y_n|<|x_1y_1|<|y_1|<\cdots< |x_ny_n|<|y_n|<1\end{align} where the trace (\ref{traceform2}) converges, the series expansion of the $\eta$-weighted $n$-point function is \begin{align}\label{firstanswer} \begin{aligned}& F_{N,\eta}\begin{pmatrix} z_1 & \dots & z_n \\ r_1 & \dots & r_n\end{pmatrix} = F_{N,\eta}(\emptyset) \prod_i \frac{1}{\vartheta_0(x_i)}\,[y_1^0\dots y_n^0]  \\ \prod_{i<j}&\, \frac{\vartheta_{r_i-r_j}(x_iy_i/x_jy_j)\vartheta_{r_i-r_j}(y_i/y_j)}{\vartheta_{r_i-r_j}(x_iy_i/y_j)\vartheta_{r_i-r_j}(y_i/x_jy_j)}\cdot \prod_{i,r} \left[\frac{\vartheta_r(x_iy_i)}{\vartheta_r(y_i)}\right]^{M(r_i-r)}.\end{aligned}\end{align} \end{theorem}

We have suppressed the second variable $q$, which appears in all of the above theta functions. We now investigate the modular properties of $\vartheta_r(x,q)$ and $F_{N,\eta}(\emptyset)$. Let $\Lambda:=2\pi i\Z\oplus 2\pi i \tau\Z$. We have the Taylor series expansions \begin{align*} -\ln \frac{\vartheta_r(e^z,q)}{z^{\delta_{r0}}} & =\sum_{k\geq 1} E_k^{r/N}(q)\frac{z^k}{k!}\hspace{10pt}\textrm{for }r\not\equiv 0\,(N) 
%\\ -\ln\frac{\vartheta_0(e^z,q)}{z} & =\sum_{k\geq 1} E_r^k(q)\frac{z^k}{k!}\hspace{10pt}.
\end{align*} where $E_k^{r/N}(q)$ is the level $\Gamma_1(N)$ Eisenstein series \begin{align*} E_k^{r/N}(q)&= (k-1)! \frac{(-1)^k}{N}\sum_{s\,(N)}\sum_{\lambda\in \Lambda} \frac{\zeta^{-rs}}{(z+\frac{2\pi si}{N})^k} \\ &=\sum_{d\in \Z\backslash \{0\}} \!\!\zeta^{-ds}\left(\frac{-N}{2\pi i d}\right)^k+\sum_{n\geq 1}q^n\!\!\!\!\sum_{\substack{ m\mid n \\ m\equiv r\,(N)}}\!\!\!\!{\rm sgn}(m)m^{k-1},\end{align*} see \cite{diashur}. Hence, the $z^k$ coefficient of $\vartheta_r(e^z,q)$ is a (quasi-)modular form of weight $k-\delta_{r0}$. By the computation from Appendix B, \begin{align*}F_{N,\eta}(\emptyset)=\frac{q^{|\eta|/N}\prod_{m\geq 1}\prod_{r,s} (1-q^{N(m-1)+[r-s]})^{M(r)M(s)}}{tr_{\mathbb{L}}\,q^{H/N}\mathfrak{W}_N \prod_m (1-q^{Nm})}.\end{align*} We may easily compute that $tr_\mathbb{L}\, q^{H/N}\mathfrak{W}_N =\prod_{m\geq 1}(1-q^m)^{-1/N}$. Separating out the terms in $M(r)$ that don't depend on $C_r$, we have $$F_{N,\eta}(\emptyset)=q^{|\eta|/N}\prod_{t\,(N)} (1-q^{N(m-1)+[t]})^{D(t)}$$ where we define the function $D\,:\,\Z/N\Z\rightarrow \Z$ by $$D(t):=C_{t-1}-C_t+C_{-t-1}-C_{-t}+\sum_{r\,(N)} (C_{r-1}-C_r)(C_{r-t-1}-C_{r-t}).$$

\begin{proposition} Let $\eta$ be an $N$-core with charges $C_r:=C_r(\eta)$. The function $F_{N,\eta}(\emptyset)$ is a modular form for $\Gamma(N)$ of weight $$C_{-1}-C_0 + \tfrac{1}{2}\sum_{r\,(N)} (C_{r-1}-C_r)^2.$$ \end{proposition}

\begin{proof} First, we observe that $D(t)=D(-t)$ and that $D(t)$ is even whenever $t\equiv -t\,(N)$. Define the function $$f_t(q):=q^{N(([t]/N)^2-[t]/N+1/6)/2 }\prod_{m\geq 1} (1-q^{N(m-1)+[t]})(1-q^{N(m-1)+[-t]}).$$ Then by Corollaries 2 and 3 of \cite{yang}, mildly generalized to allow $t\equiv 0 \,(N)$, $$\prod_{0\leq t\leq \lfloor\frac{N}{2}\rfloor}f_t(q)^{e_t}$$ is a modular function for $\Gamma(N)$ of weight $e_0$ whenever the conditions \begin{align*}\sum e_t&\equiv 0\textrm{ mod }12 \\ \sum te_t&\equiv 0\textrm{ mod }2\end{align*} are satisfied. Furthermore for odd $N$, the second condition is unnecessary. Note that the first is satisfied for any $N$-core $\eta$ because $\sum e_t=\frac{1}{2}\sum D(t)=0$. Some manipulations mod $2$ verify the second condition for $N$ even. Thus, we have that $F_{N,\eta}(\emptyset)$ is, up to some power of $q$, a modular function of weight $$ C_{-1}-C_0 + \tfrac{1}{2}\sum_{r\,(N)} (C_{r-1}-C_r)^2 .$$

We now note that the factor $q^{|\eta|/N}$ accounts exactly for the factors of $q$ appearing in $f_t(q)$. Equivalently $$|\eta| = \tfrac{N^2}{4}\sum_t D(t)\left(\tfrac{[t]^2}{N^2}-\tfrac{[t]}{N} +\tfrac{1}{6}\right).$$ This can be explicitly verified using the formula (\ref{form}), or by discrete integration. Finally, by the same argument as in Section 3.2.3 of \cite{eo2}, the growth rate of the coefficients of $F_{N,\eta}(\emptyset)$ is at most polynomial. Hence $F_{N,\eta}(\emptyset)$ is in fact of modular form for $\Gamma(N)$, rather than just a modular function. \end{proof}

Since $F_{N,\eta}(\emptyset)$ and the Taylor coefficients of $\vartheta_r(e^z,q)$ are quasimodular for $\Gamma(N)$, Theorem \ref{modularity} follows by standard coefficient extraction techniques, see e.g. Section 5 of \cite{engel} or Section 3 of \cite{eo2}. 

\begin{remark} In the author's opinion, the moral of this rather complicated computation is that the rotation operator $\mathfrak{R}_N$ acting on the tensor product $\mathbb{L}^{\otimes N}$ is the key to understanding to the elliptic orbifold of order $N$, and more generally, $N$-ic differentials $(\Sigma,\omega)$ on Riemann surfaces. Given such a differential, there is a cyclic degree $N$ cover of it $\pi\,:\,\widetilde{\Sigma}\rightarrow \Sigma$ which trivializes the $\langle \zeta_N\rangle$ monodromy of the flat metric. It feels intuitively obvious that $\mathfrak{R}_N$ is in some way related to the deck transformation of the map $\pi$, but making this connection mathematically rigorous is difficult. \end{remark}

\section*{Appendix A.}

\begin{figure}[H]
\makebox[\textwidth][c]{\includegraphics[width=5.5in]{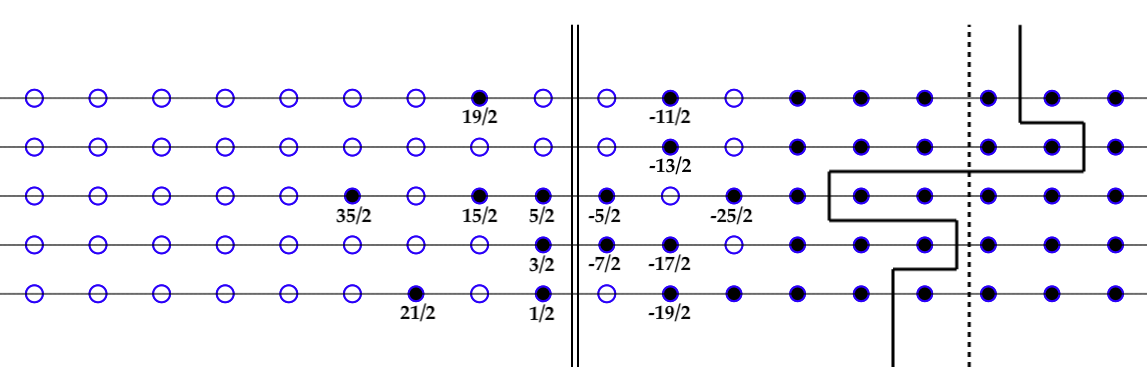}}
\caption{A plot of the values of $\lambda_i-i+\tfrac{1}{2}$ for $\lambda$ whose $r/5$-charges are $(C_0,C_1,C_2,C_3,C_4)=(1,0,2,-2,-1).$ The single solid line forms the edge of the jagged $5$-abacus; we choose it so that each string has the same number of beads to the left, and all slots from the leftmost part of the jagged edge to the right are filled.}
\label{fig2}
\end{figure}

In this appendix, we complete the proof of Theorem \ref{big} by showing that the product (\ref{superprod}) is equal to ${\bf w}_{N,\eta}(\lambda)$. Suppose that $r\equiv 1\textrm{ mod }N$. The number of times $r$ appears in (\ref{superprod}) is  \begin{align}\label{rhook}\begin{aligned} 
\#\{i\,\big{|}\,&0 \leq r+NC_{\mu_i-i-1}-NC_{\mu_i-i}\leq \mu_i-i+N\ell\} \\ & -\#\{(i,j) \,\big{|}\,r = \mu_j-\mu_i+i-j+NC_{\mu_j-j}-NC_{\mu_i-i}\}
% \\ \#\{j\,\big{|}\,&0 \leq r-NC_{[j+1]}+NC_{[j]}\leq \mu_j-j+N\ell\} \\ & -\#\{(i,j) \,\big{|}\,r = \mu_i-\mu_j+j-i+NC_{[i]}-NC_{[j]}\} &&\textrm{if }r\equiv -1\textrm{ mod }N
\end{aligned}\end{align} with a similar formula for $r\equiv -1\textrm{ mod }N$. Note that $i$ and $j$ are subject to the conditions $\mu_k-k+N\ell\geq 0$ for $k=i,j$. We visualize the above counts in terms of the {\it jagged $N$-abacus}, see Figure \ref{fig2} for a definition. Note that $v_\mu=S_\eta^{-1}v_\lambda$ results from shifting each mod $N$ substring of $v_\lambda$ so that the jagged line moves to the dotted line. Assume $r$ is not too small. Then the first term in (\ref{rhook}) can be written as \begin{align*}\#\{i\,\big{|}\,NC_{\mu_i-i-1}\leq \mu_i-i+N\ell+NC_{\mu_i-i}-r\}.\end{align*} Because $S_\eta v_\mu = v_\lambda$, we have that \begin{align*}\{\lambda_i-i\}=\{\mu_i-i+NC_{\mu_i-i}\}\end{align*} so the first term of (\ref{rhook}) exactly counts the beads $\underline{\lambda_i-i+\frac{1}{2}}$ to the left of the jagged line for which $\underline{\lambda_i-i-r+\frac{1}{2}}$ still lies on the left of the jagged line.

Now we analyze the second term in (\ref{rhook}). By the same substitution, it is \begin{align*}\#\left\{(i,j)\,\big{|}\,r=\lambda_j-\lambda_i+i-j\textrm{ with } \lambda_k-k+N\ell \geq NC_{\mu_k-k}\textrm{ for }k=i,j\right\}.\end{align*} Equivalently, it is the number of pairs of beads, both to the left of the jagged line, whose index differs by $r$. Subtracting the second term from the first therefore counts the number of beads $\underline{\lambda_i-i+\tfrac{1}{2}}$ to the left of the jagged line, such $\underline{\lambda_i-i-r+\tfrac{1}{2}}$ is also to the left of the jagged line, but is unoccupied by a bead. This is exactly the number of $r$-hooks of $\lambda$. The analysis of the $r\equiv -1\textrm{ mod }N$ terms is similar, and $r\equiv 0$ is still simpler, because all the action takes place on a single string, and the shifts play no role.

The above argument fails for small $r$ because we didn't impose \begin{align*}0\leq r + NC_{\mu_i-i}-NC_{\mu_i-i}\end{align*} in the first term of (\ref{rhook}), and have furthermore ignored the fact that $r$ can be negative and still satisfy this bound. For instance, in Figure \ref{fig1}, consider the bead $\underline{\lambda_i-i+\tfrac{1}{2}}=\underline{-\frac{5}{2}}$ and suppose $r=6$. Both $\underline{-\tfrac{5}{2}}$ and $\underline{-\tfrac{17}{2}}$ lie to the left of the jagged line, but this term is not counted in (\ref{rhook}) because $0>6+5\cdot 0-5\cdot 2=-4$ and thus in $v_\mu$ the difference in indices is negative. In fact, no beads on the $3$rd string will contribute for $r=6$. On the other hand, there is a contribution which counts towards $r=-6$ on the $4$th string coming from $\underline{\lambda_i-i+\frac{1}{2}} = \underline{-\tfrac{17}{2}}$ because when shifted to $v_\mu$, the difference in indices becomes positive.

These effects nearly cancel each other, but the $3$rd string is missing six terms for $r=6$ whereas the $4$th string only contributes an extra five $r=-6$ terms. This discrepancy is exactly the number of beads to the {\it right} of the jagged line, such that decreasing the index by $6$ lies to the {\it left} of the jagged line, corresponding in this case to the beads $\underline{-\frac{45}{2}}$ and $\underline{-\frac{57}{2}}$ on the $3$rd and $4$th strings. These discrepancies correspond exactly the $6$-hooks of $\eta=\lambda^{\textrm{mod }5}$. The same argument applies in the general case in exactly the same way. Hence by (\ref{superprod}) we have \begin{align*}\langle v_\lambda\,|\,\mathfrak{W}_{N,\eta}\,|\,v_\lambda\rangle=\frac{\displaystyle \pm\prod_{r\equiv\pm 1} r^{\#\{r\textrm{-hooks of }\lambda\}}}{\displaystyle\prod_{r\equiv \pm 1}r^{\#\{r\textrm{-hooks of }\eta\}}  \displaystyle\prod_{r\equiv 0} (r^2)^{\#\{r\textrm{-hooks of }\lambda\}}}={\bf w}_{N,\eta}(\lambda).\end{align*} \vfill

\pagebreak

\section*{Appendix B.}

In this appendix, we evaluate the operator trace (\ref{traceform2}) which gives the $\eta$-weighted $n$-point function. By Lemma \ref{rotator} this trace is:

 \begin{align*} &F_{N,\eta}(z,r) = \frac{1}{tr_\mathbb{L}\,q^{H/N}\mathfrak{W}_N \displaystyle\prod (x_i^{1/2}-x_i^{-1/2})}\prod_{m\geq 1} \frac{1}{1-q^{Nm}}\,\,[y_1^0\dots y_n^0]  \\
& \prod_{m\geq 1} \prod_{i<j}{\rm exp}\bigg(\!\!-\frac{((x_iy_i)^m-y_i^m)((x_jy_j)^{-m}-y_j^{-m})\delta_{r_i,r_j}}{m} \!\bigg)\cdot \\ 
&\prod_{m\geq 1} \prod_i {\rm exp}\bigg(\!\!-\frac{((x_iy_i)^m-y_i^m)(C_{r_i-1}-C_{r_i}+\delta_{r_i,0}-\tfrac{1}{N})}{m} \!\bigg)\cdot \\
& \prod_{m\geq 1} \prod_{i,j}{\rm exp}\bigg(\frac{-q^{[r_i-r_j]m}}{m(1-q^{Nm})}((x_iy_i)^m-y_i^m)((x_jy_j)^{-m}-y_j^{-m})\bigg) \cdot \\
& \prod_{m\geq 1} \prod_{i,r}{\rm exp}\bigg(\frac{-q^{[r_i-r_r]m}}{m(1-q^{Nm})}((x_iy_i)^m-y_i^m)(C_{r-1}-C_r+\delta_{r0}-\tfrac{1}{N})\bigg)\cdot \\
& \prod_{m\geq 1} \prod_{i,r}{\rm exp}\bigg(\frac{-q^{[r-r_i]m}}{m(1-q^{Nm})}((x_iy_i)^{-m}-y_i^{-m})(C_{r-1}-C_r+\delta_{r0}-\tfrac{1}{N})\bigg)\cdot \\
& \prod_{m\geq 1} \prod_{r,s}{\rm exp}\bigg(\frac{-q^{[s-r]m}}{m(1-q^{Nm})}(C_{r-1}-C_r+\delta_{r0}-\tfrac{1}{N})(C_{s-1}-C_s+\delta_{s0}-\tfrac{1}{N})\bigg). \end{align*} where $[r]$ denotes the unique representative of the congruence class of $r$ such that $0<[r] \leq N$. The second line of the formula comes from commuting the $\mathcal{E}_0^{r_i/N}(z_i)$ operators until they are normal-ordered. Similarly, the third line comes from commuting the raising operators in $\Psi_\eta$ past the $\mathcal{E}_0^{r_i/N}(z_i)$ operators. The fourth line comes from the trace of the normal-ordered products of $\mathcal{E}_0^{r_i/N}(z_i)$. The fifth and sixth lines come from the trace of the normal-ordered products of $\Psi_\eta$ and $\mathcal{E}_0^{r_i/N}(z_i)$ operators. Finally, the seventh line is the contribution to the trace from $\Psi_\eta$ with itself.

%Note that we have the partial fraction decomposition \begin{align}\frac{q^r}{1-q^N} = \frac{1}{N}\sum_{s\,(N)} \zeta^{-rs}\frac{\zeta_N^sq}{1-\zeta_N^sq}\end{align} for $0<s\leq N$.
To simplify this truly gigantic formula, we use the identities  \begin{align*} \prod_{m\geq 1} {\rm exp}\bigg(\pm \frac{C^m}{m} \bigg)&=(1-C)^{\mp 1} \\ \prod_{m\geq 1}  {\rm exp}\bigg(\frac{\pm C^mq^m}{m(1-q^m)} \bigg)&=\prod_{m\geq 1} (1-Cq^m)^{\mp 1}.\end{align*} Let $M(r):=C_{r-1}-C_r+\delta_{r0}-\tfrac{1}{N}$. Then the trace is

 \begin{align*} &F_{N,\eta}(z,r)= \frac{1}{tr_\mathbb{L}\,q^{H/N}\mathfrak{W}_N \displaystyle\prod (x_i^{1/2}-x_i^{-1/2})}\prod_{m\geq 1} \frac{1}{1-q^{Nm}}\,\,[y_1^0\dots y_n^0]  \\
&\prod_{\substack{ i<j \\ r_i=r_j}}\frac{(1-x_iy_i/x_jy_j)(1-y_i/y_j)}{(1-x_iy_i/y_j)(1-y_i/x_jy_j)}\cdot \\
&\prod_i \left(\frac{1-x_iy_i}{1-y_i}\right)^{M(r_i)}\cdot \\ 
& \prod_{m\geq 1} \prod_{i,j} \left[\frac{(1-(x_iy_i/x_jy_j)q^{N(m-1)+[r_i-r_j]})(1-(y_i/y_j)q^{N(m-1)+[r_i-r_j]})}{(1-(x_iy_i/y_j)q^{N(m-1)+[r_i-r_j]} )(1-(y_i/x_jy_j)q^{N(m-1)+[r_i-r_j]})}\right]\cdot \\
& \prod_{m\geq 1}\prod_{i,r}  \left[\frac{1-x_iy_iq^{N(m-1)+[r_i-r]}) }{1-y_i q^{N(m-1)+[r_i-r]}}\right]^{M(r)}\cdot \\
& \prod_{m\geq 1}\prod_{i,r} \left[\frac{1-(1/x_iy_i)q^{N(m-1)+[r-r_i]}}{1-(1/y_i) q^{N(m-1)+[r-r_i]}}\right]^{M(r)}\!\!\!\!\!\cdot \\
& \prod_{m\geq 1}\prod_{r,s}(1-q^{N(m-1)+[r-s]})^{M(r)M(s)} \end{align*}  
where each of the seven lines in the formula above correspond to the seven lines on the previous page. These factors all combine into a product of shifted theta functions in Theorem \ref{theanswer}.

%We may now introduce new variables $(w_0,\dots,w_{N-1})$ to keep track of $\eta$ in the larger generating function \begin{align}\mathbb{F}_N(z_1,\dots,z_n\,\big{|}\,r_1,\dots,r_n\,\big{|}\big{|}\,w_0,\dots,w_{N-1}):=\sum_\eta F_{N,\eta}\prod_r w_r^{C_r(\eta)}\end{align}

%\begin{remark} There are subtle sign issues which arise when comparing the action of the operators $\alpha_k^{r/N}$ on $\mathbb{L}$ and the action of $\alpha_k$ on the $r$th tensor factor of $\mathbb{L}^{\otimes N}$. To correct these issues, it is useful to put a sign into the isomorphism (\ref{tensordecomp}) as follows: \end{remark}


\begin{thebibliography}{99}

\bibitem{bo}
S. Bloch and A. Okounkov.
\emph{The character of the infinite wedge representation,}
Advances in Mathematics 149.1: 1-60,
2000.

\bibitem{diashur}
F. Diamond and J. Shurman.
\emph{A first course in modular forms,}
Vol. 140. New York: Springer,
2005.

\bibitem{eo1}
A. Eskin and A. Okounkov.
\emph{Asymptotics of numbers of branched coverings of a torus and volumes of moduli spaces of holomorphic differentials},
Inventiones Mathematicae 145.1: 59-103,
2001.

\bibitem{eo2}
A. Eskin and A. Okounkov.
\emph{Pillowcases and quasimodular forms},
Algebraic geometry and number theory,
Birkh\"auser Boston: 1-25,
2006. 

\bibitem{engel}
P. Engel.
\emph{Hurwitz Theory of Elliptic Orbifolds,}
arXiv:1706.06738,
2017.

\bibitem{engsmi}
P. Engel and P. Smillie.
\emph{The number of non-negative curvature triangulations of $S^ 2$,} Geometry \& Topology 22.5,
2018.

\bibitem{cores}
F. Garvan, D. Kim, and D. Stanton.
\emph{Cranks and t-cores},
Inventiones Mathematicae 101.1: 1-17,
1990.

\bibitem{op}
A. Okounkov and R. Pandharipande.
\emph{Gromov-Witten theory, Hurwitz theory, and completed cycles,}
Annals of mathematics 163.2: 517-560,
2006.

\bibitem{rz2}
R. Rios-Zertuche.
\emph{An introduction to the half-infinite wedge.}
Mexican Mathematicians Abroad 657: 197,
2016.

\bibitem{yang}
Y. Yang.
\emph{Transformation formulas for generalized Dedekind eta functions,}
Bulletin of the London Mathematical Society 36.5: 671-682,
2004.

\end{thebibliography}
\end{document}